\documentclass{amsart}

\usepackage{amssymb}

\theoremstyle{plain}
\newtheorem{bigthm}{Theorem}

\newtheorem{theorem}{Theorem}[section]
\newtheorem{lemma}[theorem]{Lemma}
\newtheorem{proposition}[theorem]{Proposition}
\newtheorem*{question}{Question}
\newtheorem*{theorem*}{Theorem}
\newtheorem{corollary}[theorem]{Corollary}

\theoremstyle{definition}
\newtheorem{definition}[theorem]{Definition}
\newtheorem*{example}{Example}

\theoremstyle{remark}

\newtheorem*{remarks}{Remarks}

\numberwithin{equation}{section}

\newcommand{\norm}[1]{\left\Vert #1\right\Vert}

\begin{document}

\title{Sets of k-recurrence but not (k+1)-recurrence}

\author{Nikos Frantzikinakis}
\address[Nikos  Frantzikinakis]{Department of Mathematics\\
  McAllister Building\\
  Pennsylvania State University\\
  University Park, PA \\
  16802 \\ USA} \email{nikos@math.psu.edu}

\author{Emmanuel Lesigne}
\address[Emmanuel Lesigne]{Universit\'e Fran\c{c}ois Rabelais de Tours \\
Laboratoire de Math\'ematiques et Physique Th\'eorique (UMR CNRS 6083)\\
Facult\'e des Sciences et Techniques\\
Parc de Grandmont \\
37200 Tours\\ France}
\email{lesigne@univ-tours.fr}

\author{M\'at\'e Wierdl}
\address[M\'at\'e Wierdl]{Department of Mathematical Sciences\\
  University of Memphis\\
  Memphis, TN \\ 38152 \\ USA } \email{mw@csi.hu}

\dedicatory{We dedicate this paper to Y. Katznelson.  Our work began
  at the conference organized for his 70th birthday, and we wish to
  honor him for his fundamental contribution to ergodic theory.}

\begin{abstract}
  For every $k\in \mathbb{N}$, we produce a set of integers which is
  $k$-recurrent but not $(k+1)$-recurrent.  This extends a result of
  Furstenberg who produced a $1$-recurrent set which is not
  $2$-recurrent.  We discuss a similar result for convergence of
  multiple ergodic averages.  Finally, we also point out a
  combinatorial consequence related to Szemer\' edi's theorem.
\end{abstract}

\maketitle

\tableofcontents

\section{Introduction and main results}

In his seminal paper \cite{Fu1}, Furstenberg gave an ergodic theoretic
proof of the famous theorem of Szemer\'edi claiming that every integer
subset with positive density contains arbitrarily long arithmetic
progressions. Furstenberg proved this by showing the following
multiple recurrence property for measure preserving systems:

\begin{theorem}[{\bf Furstenberg}]
  Let $(X,\mathcal{X},\mu,T)$ be a finite measure preserving system
  and $A\in \mathcal{X}$ be a set with $\mu(A)>0$. Then for every
  $k\in \mathbb{N}$, there exists $n\in\mathbb{N}$ such that $$
  \mu(A\cap T^{-n}A\cap\cdots \cap T^{-nk}A)>0.  $$
\end{theorem}

This motivated the following definition:

\begin{definition}
  Let $(X,\mathcal{X},\mu,T)$ be a probability preserving system.  We
  say that $S\subset \mathbb{N}$ is a set of \emph{$k$-recurrence for
    the transformation $T$} if for every $A\in\mathcal{X}$ with
  $\mu(A)>0$, there exists $n\in S$ such that
  \begin{equation*}
    \mu(A\cap T^{-n}A  \cap \cdots \cap T^{-kn}A)>0.
  \end{equation*}
  We say that $S\subset \mathbb{N}$ is a set of \emph{$k$-recurrence}
  if it is a set of $k$-recurrence for every probability preserving
  transformation.

\end{definition}

If $S$ is infinite then the difference set $S-S=\{s_1-s_2\colon s_1,
s_2\in S\}$ is easily shown to be a set of $1$-recurrence. By
appropriately choosing $S$, Furstenberg constructed, in \cite[pages
177-178]{Fu2}, a set of $1$-recurrence that is not a set of
$2$-recurrence.  Constructing sets of $2$-recurrence is much harder,
in fact all the examples known turned out to be sets of $k$-recurrence
for every $k$.  This raised the question, first stated explicitly by
Bergelson in \cite{B}:
\begin{question}
  Let $k\geq 2$ be an integer. Does there exist a set of
  $(k-1)$-recurrence that is not a set of $k$-recurrence?
\end{question}
The main objective of this article is to show that the answer is
affirmative. The examples that we construct are very explicit:
\begin{bigthm}\label{sec:intr-main-results}
  Let $k\geq 2$ be an integer and $\alpha\in\mathbb{R}$ be irrational.
  We define $$
  S_k=\big\{ n\in\mathbb{N} \colon
  \{n^k\alpha\}\in[1/4,3/4]\big\}, $$
  where $\{a\}$ denotes the
  fractional part of $a$. Then $S_k$ is a set of $(k-1)$-recurrence
  but not a set of $k$-recurrence.
\end{bigthm}
It will appear in the proof that not only the set $S_k$ is a set of
$(k-1)$-recurrence for powers of a single transformation, but that it
is a set of $(k-1)$-recurrence for families of commuting
transformations (see definition in Section~\ref{S:2}).

We also answer the corresponding question for sets of $k$-convergence:
\begin{definition}
  A set $S=\{a_1<a_2<\ldots\}\subset \mathbb{N}$ is called a set of
  \emph{$k$-convergence} if for for every probability preserving
  system $(X,\mathcal{X},\mu,T)$ and functions $f_1,\ldots,f_k\in
  L^\infty(\mu)$, the averages $$
  \frac{1}{N} \sum_{n=1}^N
  T^{a_n}f_1\cdot \ldots \cdot T^{ka_n}f_{k} $$
  converge in $L^2(\mu)$
  as $N\to \infty$.
\end{definition}

Host and Kra in \cite{HK} (see also Ziegler's work in \cite{ziegler}
for an alternative proof) showed that $S=\mathbb{N}$ is a set of
$k$-convergence for every $k$. We show:

\begin{bigthm}\label{sec:intr-main-results-1}
  Let $k\geq 2$ be an integer and $\alpha \in\mathbb{R}$ be
  irrational. Let $$
  I_j=
  \begin{cases}
    n\in[2^{j},2^{j+1}] \colon \{n^ka\}\in[1/10,2/10] \text{ if } j
    \text{ is even}\\
    n\in[2^{j},2^{j+1}] \colon \{n^ka\}\in[5/10,6/10] \text{ if } j
    \text{ is odd}
  \end{cases}
  $$
  and define $$
  S'_k=\bigcup_{j=1}^\infty I_j.  $$
  Then $S'_k$ is a
  set of $(k-1)$-convergence but not a set of $k$-convergence.
\end{bigthm}

It will be clear from the proof that $S'_k$ is also a set of
$(k-1)$-recurrence but not a set of $k$-recurrence.

The strategy of the proof of the theorems is as follows: In
Section~\ref{S:1} we use some elementary considerations in order to
show that $S_{k}$ is not a set of $k$-recurrence and $S'_k$ is not a
set of $k$-convergence. The basic observation is that if $S$ is a set
of $k$-recurrence/convergence then the set consisting of the $k$-th
powers of elements of $S$ has good $1$-recurrence/convergence
properties. In Section~\ref{S:2} we prove a multiple ergodic theorem
(Proposition~\ref{P:2}) that enables us to show that $S_k$ is a set of
$(k-1)$-recurrence and $S'_k$ is a set of $(k-1)$-convergence.

Finally, in Section~\ref{S:3} we derive a combinatorial consequence.

\section{Bad sets  for $k$-recurrence and
$k$-convergence}\label{S:1}
We will use the following elementary fact:
\begin{lemma}\label{L:1}
  Let $k\in\mathbb{N}$. Then there exists a nonzero integer $m$ and
  integers $l_1,\ldots, l_k$ such that
  \begin{align*}
    &l_1+2l_2+\cdots+k l_k=0\\
    & \qquad \qquad \vdots  \\
    &l_1+2^{k-1}l_2+\cdots + k^{k-1}l_k=0\\
    &l_1+2^{k}l_2+\cdots +k^kl_k = m.
  \end{align*}
\end{lemma}
\begin{proof}
  The corresponding (Vandermonde) determinant is nonzero, so for $m=1$
  we can find a rational solution for the system. After multiplying by
  an appropriate nonzero integer we get an integer solution for the
  advertised system.
\end{proof}

The following result will enable us to show that the set $S_k$ is bad
for $k$-recurrence and $S'_k$ is bad for $k$-convergence. To better
illustrate the idea, after proving the proposition, we explain how the
argument works for $k=2$.

\begin{proposition}\label{P:1}
  If $S=\{a_1<a_2<\ldots\}\subset \mathbb{N}$ is a set of
  $k$-recurrence then $S^k=\{a_1^k<a_2^k<\cdots\}$ is a set of
  $1$-recurrence for all circle rotations, and if $ S $ is a set of $
  k $-convergence then $S^k$ is a set of $1$-convergence.
\end{proposition}

\begin{remarks}
  Note that in the above proposition we only claim that $S^k$ is a set
  of $1$-recurrence for rotations of the circle. It is clear that the
  argument we give in the proof below can be extended to show that if
  $S$ is a set of $k$-recurrence, then $S^k$ is a set of recurrence
  for all translations of multidimensional tori.  In fact, it is an
  \emph{unsolved problem} whether a set $S$ being of $k$-recurrence
  implies that $S^k$ is a set of $1$-recurrence.

  Here is a related unsolved problem of Katznelson from \cite{K}: is
  it true that a set of recurrence for all translations of
  multidimensional tori is, in fact, a set of
  \emph{topological} recurrence?
\end{remarks}

\begin{proof}[Proof of Proposition~\ref{P:1}]
  (i) Let $S$ be a set of $k$-recurrence. It suffices to show that for
  every $\alpha\in [0,1)$ and $\varepsilon>0$ there exists
  $n\in\mathbb{N}$ such that $\norm{a_n^k \alpha}\leq \varepsilon$,
  where $\norm{a}$ is the distance of $a$ to the closest integer. To
  do this we use the $k$-recurrence property for some appropriately
  chosen system. We define the measure preserving transformation $R$
  acting on $\mathbb{T}^k$ with the Haar measure $\lambda$ as follows:
  the $i$-th coordinate of $R(t_1,\ldots,t_k)$ is $$
  t_i+\binom{i}{1}t_{i-1}+\binom{i}{2}t_{i-2}+\cdots+\binom{i}{i-1}t_1+\alpha.
  $$
  The last coordinate of $R^{jn}(t_1,\ldots,t_k)$ is
  \begin{equation}\label{E:c}
    c_{j,n}=t_k+\binom{k}{1}(jn)t_{k-1}+\binom{k}{2}(jn)^2t_{k-2}+\cdots+\binom{k}{k-1}(jn)^{k-1}t_1+
    (jn)^k\alpha.
  \end{equation}
  By Lemma~\ref{L:1} there exist $l_1,\ldots,l_k\in \mathbb{Z}$ such
  that
  \begin{equation}\label{E:1}
    l_1 c_{1,n}+\cdots +l_kc_{k,n}=m n^k \alpha+(l_1+\cdots +l_k)t_k
  \end{equation}
  holds for some nonzero integer $m$, and for all $n\in\mathbb{N}$,
  $\alpha\in [0,1)$. Let $U_\varepsilon=B(0,\varepsilon/(2M))$ where
  $M=|l_1|+\cdots+|l_k|$. If $S$ is a set of $k$-recurrence then there
  exists an $n_0\in S$ such that $$
  U_\varepsilon\cap
  R^{-n_0}U_\varepsilon \cap \cdots \cap R^{-kn_0}U_\varepsilon \neq
  \varnothing.  $$
  Let $(t_1,\ldots,t_k)$ be an element of the
  intersection and $c_{j,n}$, $j=1,\ldots,k$, be given by \eqref{E:c}.
  We have $\norm{c_{j,n_0}}\leq \varepsilon/2M$, $j=1,\ldots, k$, and
  $\norm{t_k}\leq \varepsilon/(2M)$. Using \eqref{E:1} we get that
  $\norm{n_0^k (m \alpha)}\leq \varepsilon$. Since $m$ does not depend
  on the choice of $\alpha$, and for every nonzero integer $m$ the map
  $\alpha \to m\alpha$ is onto, the result follows.

  (ii) Suppose that $S=\{a_1<a_2<\cdots\}$ is a set of
  $k$-convergence. Let $l_1,\ldots, l_k, m$ be as in the proof part
  (i), and let $R$ be the transformation defined there. If
  $f_j(t_1,\ldots,t_k)=e(l_j t_k)$, $j=1,\ldots, k$, then the averages
  $$
  \frac{1}{N} \sum_{n=1}^N R^{a_n}f_1\cdot \ldots \cdot
  R^{ka_n}f_{k}= e\big((l_1+\cdots+l_k)t_k\big) \cdot \frac{1}{N}
  \sum_{n=1}^N e(a_n^km\alpha) $$
  converge in $L^2(\lambda)$ as
  $N\to\infty$. Hence, for every $\alpha\in [0,1)$ the averages $$
  \frac{1}{N} \sum_{n=1}^N e(a_n^k\alpha) $$
  converge as $N\to\infty$.
  The spectral theorem gives that for every measure preserving system
  $(X,\mathcal{X},\mu,T)$ and $f\in L^2(\mu)$ the averages $$
  \frac{1}{N} \sum_{n=1}^N T^{a_n^k}f $$
  converge in $L^2(\mu)$ as
  $N\to\infty$, completing the proof.
\end{proof}

We now give an example illustrating how the argument of
Proposition~\ref{P:1} works for $k=2$:

\begin{example}
  (i) Suppose that $S$ is a set of $2$-recurrence. Let $\alpha\in
  [0,1)$ and $\varepsilon>0$. The transformation $R\colon
  \mathbb{T}^2\to\mathbb{T}^2$ is defined by $$
  R(t_1,t_2)=(t_1+\alpha,t_2+2t_1+\alpha).  $$
  Then $$
  R^n(t_1,t_2)=(t_1+n\alpha,t_2+2nt_1+n^2\alpha).  $$
  Since $S$ is a
  set of $2$-recurrence there exists an $n_0\in S$ such that $$
  U_\varepsilon\cap R^{-n_0}U_\varepsilon \cap R^{-2n_0}U_\varepsilon
  \neq \varnothing, $$
  where $U_\varepsilon=B(0,\varepsilon/4)$.  If
  $(t_1,t_2)$ is an element of the intersection then

\begin{equation}\label{E:4}
  \norm{t_2}\leq \varepsilon/4, \ \norm{c_1}\leq \varepsilon/4, \ \
  \norm{c_2}\leq \varepsilon/4,
\end{equation}
where $$
c_1=t_2+2n_0t_1+n_0^2\alpha, \ c_2=t_2+4n_0t_1+4n_0^2\alpha.
$$
Since $$
c_2-2c_1=2n_0^2\alpha+t_2 $$
we get from \eqref{E:4} that
$\norm{n_0^2(2\alpha)}\leq \varepsilon$. This implies that $S^2$ is a
set of recurrence for circle rotations.

(ii) Suppose that $S=\{a_1<a_2<\cdots\}$ is a set of $2$-convergence.
We define the transformation $R\colon\mathbb{T}^2\to\mathbb{T}^2$ as
in part $(i)$ and the functions $$
f(t_1,t_2)=e(-2t_2), \
f_2(t_1,t_2)=e(t_2).  $$
Then $$
\frac{1}{N}\sum_{n=1}^N
R^{a_n}f_1\cdot R^{2 a_n} f_2= e(-t_2)\cdot \frac{1}{N} \sum_{n=1}^N
e(a_n^2\alpha).  $$
Since $S$ is a set of $2$-convergence it follows
that the averages $$
\frac{1}{N}\sum_{n=1}^N e(a_n^2\alpha) $$
converge as $N\to\infty$ for every $\alpha\in [0,1)$. The spectral
theorem gives that $S^2$ is a set of $2$-convergence.
\end{example}

\begin{corollary}
  The set $S_k$ of Theorem~\ref{sec:intr-main-results} is {\bf not} a
  set of $k$-recurrence and the set $S'_k$ of
  Theorem~\ref{sec:intr-main-results-1} is {\bf not} a set of
  $k$-convergence.
\end{corollary}
\begin{proof}
  By definition, $S_k^{k}$ is not a set of recurrence for the rotation
  by $\alpha$, so by Proposition~\ref{P:1} we have that $S_k$ is not a
  set of $k$-recurrence.

  Let $S'_k=\{a_1<a_2<\cdots\}$ and $$
  A_N=\frac{1}{N}\sum_{n=1}^N
  e(a_n^{k} \alpha), \quad B_N=\frac{1}{N}\sum_{n=N+1}^{2N} e(a_n^{k}
  \alpha).  $$
  By the definition of $S'_k$ we have that for $j$ even
  the real part of $B_{2^{j}}$ is positive and bounded away from zero,
  and for $j$ odd the real part of $B_{2^{j}}$ is negative. Hence,
  the sequence $B_N$ does not converge as $N\to\infty$. Since
  $B_N=2A_{2N}-A_N$ it follows that the sequence $A_N$ does not
  converge as $N\to\infty$. By Proposition~\ref{P:1}, $S'_k$ is not a
  set of $k$-convergence.
\end{proof}

\section{Good sets  for $k$-recurrence and
$k$-convergence}\label{S:2} We will use the following elementary
lemma (\cite{B1}):
\begin{lemma}[\bf{Van der Corput}]
  Let $\{a_n\}_{n\in\mathbb{N}}$ be a bounded sequence of vectors in a
  Hilbert space. For each $m$ we set $$
  b_m=\limsup_{N-M\to\infty}\Big|\frac{1}{N-M}\sum_{n=M+1}^{N}
  <a_{n+m},a_n>\Big|.  $$
  If $$
  \limsup_{N-M\to\infty}\frac{1}{N-M}\sum_{m=M+1}^{N} b_m=0 $$
  then $$
  \lim_{N-M\to\infty}\frac{1}{N-M}\sum_{n=M+1}^{N} a_n=0 $$
  in norm.
\end{lemma}
\begin{proposition}\label{P:2}
  Let $k\in\mathbb{N}$, $T_1,\ldots,T_{k-1}$ be commuting measure
  preserving transformations acting on the probability space
  $(X,\mathcal{X},\mu)$, and $p\in\mathbb{R}[t]$ be a polynomial of
  degree $\geq k$ with irrational leading coefficient. If $g \colon
  \mathbb{T}\to \mathbb{C}$ is Riemann integrable and $f_1, \ldots,
  f_{k-1}\in L^\infty(\mu)$, then the difference $$
  \frac{1}{N-M}
  \sum_{n=M+1}^{N} T_1^nf_1\cdot \ldots \cdot T_{k-1}^nf_{k-1}\cdot
  g(p(n))- \frac{1}{N-M} \sum_{n=M+1}^{N} T_1^nf_1\cdot \ldots \cdot
  T_{k-1}^nf_{k-1} \cdot \int_{\mathbb{T}} g \ d\lambda $$
  converges
  to $0$ in $L^2(\mu)$ as $N-M\to\infty$.
\end{proposition}
\begin{proof}
  Using the standard estimation by continuous functions from above and
  below, it suffices to check the result when $ g $ is a continuous
  function. By Weierstrass approximation theorem of continuous
  functions by trigonometric polynomials, and using linearity, it
  suffices to check the result when $g(t)= e^{2\pi i l t}$ for some
  $l\in\mathbb{Z}$. The case $l=0$ is trivial. If $l\neq 0$ the
  polynomial $q(n)=lp(n)$ satisfies the assumptions of our theorem, so
  it suffices to verify the result when $g(t)=e(t)$, where
  $e(t)=e^{2\pi i t}$.

  We proceed by induction on the number of functions $k$. If $k=1$
  (empty product is $1$) then $$
  \lim_{N-M\to\infty}\frac{1}{N-M}\sum_{n=M+1}^{N}e(p(n))=0 $$
  follows
  from Weyl's uniform distribution theorem. Assume that the result is
  true for $k-1$ functions. We will verify the result for $k$
  functions. We can assume that $\norm{f_k}_\infty\leq 1$. We apply
  Van der Corput's Lemma on the Hilbert space $L^2(\mu)$ for the
  sequence of functions $$
  a_n(x)=f_1(T_1^nx)\cdot \ldots \cdot
  f_{k}(T_{k}^n x) \cdot e(p(n)).  $$
  It suffices to verify that
  for every $m\in\mathbb{N}$ the averages
  \begin{gather*}
    \frac{1}{N-M}\sum_{n=M+1}^{N} \int a_{n+m}(x)\cdot
    \overline{a_n(x)} \
    d\mu=\\
    \frac{1}{N-M} \sum_{n=M+1}^{N} \int T_1^n(T_1^mf_1 \cdot
    \bar{f_1})\cdot \ldots \cdot T_{k}^n(T_k^mf_k\cdot \bar{f_k})
    \cdot e\big(p(n+m)-p(n)\big)\ d \mu
  \end{gather*}
  go to 0 as $N-M\to \infty$.  Introducing the notation $S_i=T_i
  T_k^{-1}$, $i=1,\ldots,k-1$, and using Cauchy's inequality, it
  suffices to prove that
  \begin{equation*}
       \frac{1}{N-M} \sum_{n=M+1}^{N} S_1^n(T_1^mf_1 \cdot
      \bar{f_1})\cdot \ldots \cdot S_{k-1}^n(T_{k-1}^mf_{k-1}\cdot
      \bar{f}_{k-1}) \cdot e\big(p(n+m)-p(n)\big)
  \end{equation*}
  converge to zero in $L^2(\mu)$ as $N-M\to \infty$.  But this follows
  from the induction hypothesis since the transformations $S_i$
  commute, and for $m\in\mathbb{N}$ the polynomial $q(n)=p(n+m)-p(n)$
  has degree $\geq k-1$ and irrational leading coefficient.
\end{proof}

We remark that the non-uniform version ($M=0$) of the previous result
suffices for the proof of the next corollary\footnote{See our note at
  \url{http://www.csi.hu/mw/general_dyadic_construction_short.pdf}.}
The uniform version is only used to simplify the proof.

\begin{definition}
  We say that $S\subset \mathbb{N}$ is a set of \emph{$k$-recurrence
    for commuting transformations} if whenever $T_1,\ldots,T_k$ are
  commuting measure preserving transformations acting on the
  probability space $(X,\mathcal{X},\mu)$ and $A\in\mathcal{X}$ with
  $\mu(A)>0$, there exists $n\in S$ such that $$
  \mu(A\cap T_1^{-n}A
  \cap \cdots \cap T_k^{-n}A)>0.  $$
\end{definition}

\begin{corollary}
  The set $S_k$ of Theorem~\ref{sec:intr-main-results} is a set of
  $(k-1)$-recurrence for commuting transformations and the set $S'_k$
  of Theorem~\ref{sec:intr-main-results-1} is a set of
  $(k-1)$-convergence.
\end{corollary}
\begin{proof}
  To show that $S_k$ is a set of $(k-1)$-recurrence for commuting
  transformations we apply Proposition~\ref{P:2} for $g(t)={\bf
    1}_{[1/4,3/4]}(t)$ and $p(n)=n^k\alpha$. We get that if
  $T_1,\ldots,T_{k-1}$ are commuting measure preserving
  transformations acting on the probability space
  $(X,\mathcal{X},\mu)$, and $A\in\mathcal{B}$ with $\mu(A)>0$, then
  \begin{align*}
    &\limsup_{N\to\infty}\frac{1}{N} \sum_{n=1}^N {\bf 1}_{S_k}(n)
    \cdot\mu(A\cap T_1^{-n}A \cap \cdots \cap T_{k-1}^{-n}A)=\\
    \frac{1}{2}\cdot &\limsup_{N\to\infty}\frac{1}{N} \sum_{n=1}^N
    \mu(A\cap T_1^{-n}A \cap \cdots \cap T_{k-1}^{-n}A)>0,
  \end{align*}
  where positiveness follows from the multiple recurrence theorem of
  Furstenberg and Katznelson \cite{FK1}. Hence, $S_k$ is a set of
  $(k-1)$-recurrence.

  To show that $S'_k$ is a set of $(k-1)$-convergence we apply
  Proposition~\ref{P:2} for $T_i=T^i$, $i=1,\ldots,k-1$,
  $p(n)=n^k\alpha$, and $g={\bf 1}_{[1/10,2/10]}$ on intervals of the
  form $[2^j,2^{j}+N)$, for large $N<2^j$,when $j$ is even, and
  $g={\bf 1}_{[5/10,6/10]}$ on intervals of the form $[2^j,2^{j}+N)$,
  for large $N<2^j$, when $j$ is odd. We get that for every
  $f_1,\ldots, f_{k-1}\in L^{\infty}(\mu)$ the difference $$
  \frac{1}{N}\sum_{n=1}^N {\bf 1}_{S'_k}(n) \cdot T^nf_1\cdot\ldots
  \cdot T^{(k-1)n}f_{k-1}- \frac{1}{10}\cdot \frac{1}{N}\sum_{n=1}^N
  T^nf_1\cdot\ldots \cdot T^{(k-1)n}f_{k-1} $$
  converges to zero in
  $L^2(\mu)$ as $N\to\infty$. We know from \cite{HK} that the averages
  $$
  \frac{1}{N}\sum_{n=1}^N T^nf_1\cdot\ldots \cdot T^{(k-1)n}f_{k-1}
  $$
  converge in $L^2(\mu)$ as $N\to\infty$, so the set $S'_k$ is a
  set of $(k-1)$-convergence. This completes the proof.
\end{proof}
The reason we cannot prove that $S'_k$ is a set of $(k-1)$-convergence
for commuting transformations is that we do not yet know the analogous
convergence result for the averages $$
\frac{1}{N}\sum_{n=1}^N
T_1^nf_1\cdot\ldots \cdot T_k^{n}f_{k}.$$

\section{Combinatorial consequence}\label{S:3}

A set $S\subset\mathbb{N}$ is called intersective if for every integer
subset $\Lambda$ with positive density we have
$\Lambda\cap(\Lambda+n)\neq\varnothing$ for some $n\in S$. More
generally we define:

\begin{definition}
  A set $S\subset \mathbb{N}$ is \emph{$k$-intersective} if every
  integer subset with positive density contains at least one
  arithmetic progression of length $k+1$ and common difference in $S$.
\end{definition}

In \cite[pages 528-529]{FKO} it is shown that:

\begin{proposition}
  A set $S\subset \mathbb{N}$ is $k$-intersective if and only if it is
  a set of $k$-recurrence.
\end{proposition}

We conclude from Theorem~\ref{sec:intr-main-results} that:

\begin{corollary}
  Let $k\geq 2$. There exists a set that is $(k-1)$-intersective but
  not $k$-intersective.
\end{corollary}
\bibliographystyle{plain}
\bibliography{recurrence_submit}

\end{document}